\documentclass[12pt,a4paper]{amsart}

\usepackage{amsthm,amsmath,amssymb}
\RequirePackage[numbers]{natbib}
\RequirePackage[colorlinks,citecolor=blue,urlcolor=blue]{hyperref}
\def\bR{\mathbb R}
\def\bE{\mathbb E}

\newtheorem{Theorem}{Theorem}

\newtheorem{lemma}[Theorem]{Lemma}
\newtheorem{Corollary}[Theorem]{Corollary}
\newtheorem{proposition}[Theorem]{Proposition}
\newtheorem{Assumption}[Theorem]{Assumption}

\def\argmin{\mathop{\rm arg\,min}}

\def\rank{{\rm rank}}

\def\Sum{\overset{n}{\underset{i=1}{\sum}}}

\def\matrix{\mathbb{R}^{m_{1}\times m_{2}}}

\begin{document}

% "Title of the Paper"
\title[Rank penalized estimators]{Rank penalized estimators for high-dimensional matrices}
%\thankstext{t1}{This is an original survey paper}

% indicate corresponding author with \corref{}
% \author{\fnms{John} \snm{Smith}\thanksref{t2}\corref{}\ead[label=e1]{smith@foo.com}\ead[label=e2,url]{www.foo.com}}
% \thankstext{t2}{Thanks to somebody} 
% \address{line 1\\ line 2\\ \printead{e1}\\ \printead{e2}}
\author{Olga Klopp}
\address{Laboratoire de Statistique, CREST and University Paris Dauphine\\CREST 3, Av. Pierre Larousse 92240 Malakoff France}
\email{olga.klopp@ensae.fr}
%\and
%\author{\fnms{???} \snm{???}\ead[label=e2]{???}}
%\address{\printead{e2}}
\begin{abstract}
In this paper we consider the trace regression model. Assume that we observe  a small set of entries or linear combinations of entries of an unknown matrix $A_0$ corrupted by noise. We propose a new rank penalized estimator of $A_0$. For this estimator we establish general oracle inequality for the prediction error both in probability and in expectation. We also prove upper bounds for the rank of our estimator. Then, we apply our general results to the problems of matrix completion and matrix regression. In these cases our estimator has a particularly simple form: it is obtained by hard thresholding of the singular values of a matrix constructed from the observations.
\end{abstract}
%%\kwd{}
%
%
%% history:
%% \received{\smonth{1} \syear{0000}}
%
%%\tableofcontents
%
%%\usepackage[utf8x]{inputenc}
%%\usepackage{ucs}
%%\usepackage{amsmath}
%%\usepackage{amsfonts}
%%\usepackage{amssymb}
%%\usepackage[hypertex]{hyperref}
%
%% à utiliser pour pdf
%%
%%\usepackage[dvips,bookmarks,colorlinks=true,linkcolor=blue,urlcolor=red]{hyperref}
%
%
%% 
%%
%%\email{\href{mailto:olga.klopp@ensae.fr}{olga.klopp@math.ensae.fr}}
%
\maketitle

%\begin{keyword}[class=AMS]
%\kwd{62J99, 62H12, 60B20, 60G15}
%%\kwd{}
%%\kwd[; secondary ]{}
%\end{keyword}
%
%\begin{keyword}
%\kwd{matrix completion, low rank matrix estimation}
%\kwd{ recovery of the rank, statistical learning}
%\end{keyword}

% history:
% \received{\smonth{1} \syear{0000}}

%\tableofcontents

%\end{frontmatter}

%\usepackage[utf8x]{inputenc}
%\usepackage{ucs}
%\usepackage{amsmath}
%\usepackage{amsfonts}
%\usepackage{amssymb}
%\usepackage[hypertex]{hyperref}

% à utiliser pour pdf
%
%\usepackage[dvips,bookmarks,colorlinks=true,linkcolor=blue,urlcolor=red]{hyperref}

% 
%
%\email{\href{mailto:olga.klopp@ensae.fr}{olga.klopp@math.ensae.fr}}

%\begin{document}
%\maketitle
\section{Introduction}
In this paper we consider  the trace regression problem. Assume that we observe $n$ independent random pairs $(X_i,Y_i)$, $i=1,\dots,n$. Here  $X_{i}$ are random matrices of dimension $m_{1}\times m_{2}$ and  known distribution $\Pi_i$, $Y_i$ are random variables in $\mathbb{R}$ which satisfy
\begin{equation}
\bE \left (Y_i|X_i\right )=\mathrm{tr}(X_{i}^{T}A_{0}), \:i=1,\dots,n,
\end{equation}
where $A_{0}\in \mathbb{R}^{m_{1}\times m_{2}}$ is an unknown matrix, $\bE \left (Y_i|X_i\right )$ is the conditional expectation of $Y_i$ given $X_i$ and $\mathrm{tr}(A)$ denotes the trace of the matrix $A$. We consider the problem of estimating of $A_0$ based on the observations $(X_i,Y_i), i=1,\dots,n$.
Though the results of this paper are obtained for general $n,m_1,m_2$, our main motivation is the high-dimensional case, which corresponds to $m_1m_2\gg n$, with low rank matrices $A_0$.

Setting $\xi_{i}=Y_i-\bE \left (Y_i|X_i\right )$ we can equivalently write our model in the form
\begin{equation}\label{model}
Y_{i}=\mathrm{tr}(X_{i}^{T}A_{0})+\xi_{i}, \:i=1,\dots,n.
\end{equation}
The noise variables $(\xi_{i})_{i=1,\dots, n}$ are independent and have mean zero.

The problem of estimating low rank matrices recently generated a considerable number of works. The most popular methods are based on minimization of the  empirical risk penalized by the  nuclear norm  with  various modifications, see, for example,  \cite{Argyriou:2008:CMF:1455903.1455908,Argyriou:2010:SL:1756006.1756037,DBLP:conf/nips/ArgyriouMPY07,
MR2417263,DBLP:journals/corr/abs-0903-3131, DBLP:journals/corr/abs-1001-0339,2010arXiv1008.4886G,2009arXiv0912.5100N,2009arXiv0912.5338R}.

In this paper we propose a new estimator of $A_0$. In our construction we combine the penalization by the rank with the use of the knowledge of the distribution   $\Pi=\dfrac{1}{n}\Sum \Pi_i$. An important feature of our estimator is that in a number of interesting examples we can write it out explicitly. %It reduces to hard thresholding of a matrix constructed from the observed data $(X_i,Y_i)$.

 Penalization by the rank was previously considered in \cite{2010arXiv1004.2995B, 2010arXiv1009.5165G} for the multivariate response regression model. The criterion introduced by Bunea, She and Wegkamp in \cite{2010arXiv1004.2995B}, the rank selection criterion (RSC), minimizes the Frobenius norm of the fit plus a regularization term proportional to the rank. The rank of the RSC estimator gives a consistent estimation of the number of the singular values of the signal
 $XA_0$ above a certain noise level. Here $X$ is the matrix of predictors. %This level is given by the expectation of the largest singular value of the projection of the noise onto the column space of the predictors $X_i$.
  In \cite{2010arXiv1004.2995B} the authors also establish oracle inequalities on the mean squared errors of RSC. The paper \cite{2010arXiv1009.5165G} is mainly focused on the case of unknown variance of the noise. The author  gives a minimal sublinear penalty for RSC and provides oracle inequalities on the mean squared risks. 

The idea to incorporate the knowledge of the distribution $\Pi$ in the construction of the estimator was first introduced in \cite{Koltchinskii-Tsybakov} but with a different penalization term, proportional to the nuclear norm. In \cite{Koltchinskii-Tsybakov} the authors establish general sharp oracle inequalities for trace regression model and apply them to the noisy matrix completion problem. They also provide lower bounds.

 In the present work we consider a more general model than the model of \cite{2010arXiv1004.2995B, 2010arXiv1009.5165G}. It contains as  particular cases a number of interesting problems such as matrix completion, multi-task learning, linear regression model, matrix res\-ponse regression model. The analysis of our model requires different techniques and uses the matrix version of Bernstein's inequality for the estimation of the stochastic term, similarly to \cite{Koltchinskii-Tsybakov}. However, we use a different penalization term   than in \cite{Koltchinskii-Tsybakov} and the main scheme of our proof is quite different. In particular, we obtain a bound for the rank of our estimator in a very general setting (Theorem \ref{thm1}, (i)) and estimations for the prediction error in expectation (Theorem \ref{thm2}). Such bounds  are not available for nuclear norm penalization used in \cite{Koltchinskii-Tsybakov}. Note, however, that under very specific assumptions on $X_i$, \cite{MR2417263} shows that the rank of $A_0$ can be reconstructed exactly, with high probability, when the dimension of the problem is smaller then the sample size.

The paper is organized as follows. In Section \ref{def} we define the main objects of our study, in particular, our estimator. We also show how some well-known problems (matrix completion, column masks,''complete'' subgaussian design) are related to our model. In Section \ref{general}, we show that the rank of our estimator is bounded from above by the rank of the unknown matrix $A_0$ with a constant close to 1. In the same section we prove general oracle inequalities for the prediction error both in probability and in expectation. Then, in Section \ref{matrix} we apply these general results to the noisy matrix completion problem. In this case our estimator has a particularly simple form: it is obtained by hard thresholding of the singular values of a matrix constructed from the observations $(X_i,Y_i), i=1,\dots,n$. Moreover, up to a logarithmic factor, the rates attained by our estimator are optimal under the Frobenius risk for a simple class of matrices $\mathcal{A}(r,a)$ defined as follows: for any $A_0\in \mathcal{A}(r,a)$ the rank of $A_0$ is supposed not to be larger than a given $r$ and all the entries of $A_0$ are supposed to be bounded in absolute value by a constant $a$. Finally, in Section \ref{matrix regression}, we consider the matrix regression model and compare our bounds to those obtained in \cite{2010arXiv1004.2995B}.

\section{Definitions and assumptions}\label{def}
For $0<q\leq\infty$ the Schatten-q (quasi-)norm of the matrix $A$  is defined by
\begin{equation*}
\parallel A\parallel_q=\left (\underset{j=1}{\overset{\min(m_1,m_2)}{\Sigma}}\sigma_j(A)^{q}\right )^{1/q}\text{for}\; 0<q<\infty\;\; \text{and}\; \parallel A\parallel_\infty=\sigma_1(A),
\end{equation*}
where  $(\sigma_j(A))_j$ are the singular values of $A$ ordered decreasingly.

For any matrices $A,B\in \mathbb{R}^{m_{1}\times m_{2}}$, we define the scalar product
\begin{equation*}
\langle A,B\rangle =\mathrm{tr}(A^{T}B)
\end{equation*}
and the bilinear symmetric form
\begin{equation}\label{produit-scal}
\langle A,B\rangle_{L_{2}(\Pi)} =\dfrac{1}{n}\Sum\bE\left(\langle A,X_i\rangle\langle B,X_i\rangle\right), \quad\text{where}\quad\Pi=\dfrac{1}{n}\Sum \Pi_i.
\end{equation}
We introduce the following assumption on the distribution of the matrix $X_{i}$:
\begin{Assumption}\label{ass1} There exists a constant $\mu>0$ such that, for all matrices $A\in \mathbb{R}^{m_{1}\times m_{2}}$
$$\parallel A \parallel _{L_{2}(\Pi)}^{2}\geq \mu^{-2}\parallel A\parallel _{2}^{2}.$$
\end{Assumption}
Under Assumption 1 the bilinear form defined by \eqref{produit-scal} is a scalar product. This assumption is satisfied, often with equality, in several interesting examples such as matrix completion, column masks, ``complete'' subgaussian design.

The trace regression model is quite a general model which contains as particular cases a number of interesting problems: 
\begin{itemize}
\item
\textbf{\underline{Matrix Completion}} Assume that the design matrices $X_i$ are i.i.d uniformly distributed on the set
\begin{equation}\label{basisUSR}
\mathcal{X} = \left \{e_j(m_1)e_k^{T}(m_2),1\leq j\leq m_1, 1\leq k\leq m_2\right \},
\end{equation}
where $e_l(m)$ are the canonical basis vectors in $\bR^{m}$. Then, the problem of estimating $A_0$ coincides with the problem of matrix completion under uniform sampling at random (USR). The latter problem was studied in \cite{DBLP:journals/corr/abs-0910-1879, DBLP:journals/corr/abs-0910-0651} in the non-noisy case ($\xi_i=0$) and in \cite{2009arXiv0912.5338R, 2010arXiv1008.4886G,Koltchinskii-Tsybakov} in the noisy case. In a slightly different setting the problem of matrix completion was considered, for example, in \cite{DBLP:journals/corr/abs-1001-0339,DBLP:journals/corr/abs-0903-3131, DBLP:journals/corr/abs-0805-4471,DBLP:journals/corr/abs-0910-1879,DBLP:journals/corr/abs-0906-2027}.  

For such $X_i$, we have the relation
\begin{equation}
m_1m_2\parallel A\parallel_{L_{2}(\Pi)}^{2}=\parallel A \parallel^{2}_{2},
\end{equation}
for all matrices $A\in\matrix$.

\item\textbf{\underline{Column masks}} Assume that the design matrices $X_i$ are i.i.d. replications of a random matrix $X$, which has only one nonzero column. If the distribution of $X$ is such that all the columns have the same probability to be non-zero and the non-zero column $X_j$ is such that $\bE \left (X_jX_j^{T}\right )$ is the identity matrix, then the Assumption \ref{ass1} is satisfied with $\mu=\sqrt{m_2}$.

\item\textbf{\underline{``Complete'' subgaussian design}}  Suppose  that the design matrices $X_i$ are i.i.d. replications of a random matrix $X$ and the entries of $X$ are either i.i.d. standard Gaussian or Rademacher random variables. In both cases, Assumption \ref{ass1} is satisfied with $\mu=1$.

\item\textbf{\underline{Matrix regression}}  The matrix regression model is given by 
\begin{equation}\label{regression}
U_i=V_i\,A_0+E_i\qquad i=1,\dots, l,
\end{equation}
where $U_i$ are $1\times m_2$ vectors of response variables, $V_i$ are $1\times m_1$ vectors of predictors, $A_0$ is an unknown $m_1\times m_2$ matrix of regression coefficients and $E_i$ are random $1\times m_2$ vectors of noise with independent entries and mean zero.

We can equivalently write this model as a trace regression model. Let $U_i=(U_{ik})_{k=1,\dots, m_2}$, $E_i=(E_{ik})_{k=1,\dots, m_2}$ and $Z^{T}_{ik}=e_k(m_2)\,V_i$, where $e_k(m_2)$ are the $m_2\times 1$ vectors of the canonical basis of $\mathbb {R}^{m_2}$. Then we can write \eqref{regression} as
\begin{equation*}
U_{ik}=\mathrm{tr}(Z_{ik}^{T}A_0)+E_{ik}\qquad i=1,\dots, l \quad\text{and}\quad k=1,\dots, m_2.
\end{equation*}
Set 
%$n=lm_2$ and 
$V=\left  (V_1^{T},\dots, V_l^{T}\right  )^{T}$ and $U=\left  (U_1^{T},\dots, U_l^{T}\right  )^{T}$. Then
\begin{equation*}
\parallel A \parallel _{L_{2}(\Pi)}^{2}=\dfrac{1}{l\,m_2}\bE\left (\parallel VA\parallel^{2}_2\right ).
\end{equation*}
Assumption \ref{ass1}, which is a condition of isometry in expectation, is used in the case of random $X_i$. In the case of matrix regression with deterministic $V_i$ we do not need it, see section \ref{matrix regression} for more details.

\item\textbf{\underline{Linear regression with vector parameter}}
Let $m_1=m_2$ and $\mathbb{D}$ denotes the set of diagonal matrices of size $m_1\times m_1$. If $A$ and $X_i\in \mathbb{D}$ then the trace regression model  becomes the linear regression model with vector parameter.
\end{itemize}
%The general oracle inequalities that we prove in Section \ref{general} can be applied to these examples. 
The main motivation of this paper is the matrix completion and matrix regression problems, which we treat in Section \ref{matrix} and Section \ref{matrix regression}.

We define the following estimator of $A_0$:
\begin{equation}\label{estimator}
\hat{A}=\underset{A\in \mathbb{R}^{m_{1}\times m_{2}}}{\argmin}\left \{\parallel A\parallel _{L_{2}(\Pi)}^{2}-\Big \langle \dfrac{2}{n}\Sum Y_iX_i,A\Big \rangle+\lambda \rank (A)\right \},
\end{equation}
where $\lambda>0$ is a regularization parameter and $\rank (A)$ is the rank of the matrix $A$. 

For matrix regression problem and deterministic $X_i$, our estimator coincides with the RSC estimator:
\begin{equation*} %\label{estimator}
\begin{split}
\hat{A}&=\underset{A\in \mathbb{R}^{m_{1}\times m_{2}}}{\argmin}\left \{\dfrac{1}{l\,m_2}\parallel VA\parallel _2^{2}-\dfrac{2}{l\,m_2}\Big \langle V^{T}U,A\Big \rangle+\lambda \rank (A)\right \}\\&=
\underset{A\in \mathbb{R}^{m_{1}\times m_{2}}}{\argmin}\left \{\parallel U-VA\parallel _2^{2}+lm_2\lambda \rank (A)\right \}.
\end{split}
\end{equation*}
This estimator, called the RSC estimator, can be computed efficiently using the procedure described in \cite{2010arXiv1004.2995B}.

Under Assumption 1, the functional 
$$A\mapsto \psi(A)= \parallel A\parallel _{L_{2}(\Pi)}^{2}-\Big \langle \dfrac{2}{n}\Sum Y_iX_i,A\Big \rangle+\lambda \rank (A)$$ 
tends to $+\infty$ when $\parallel A\parallel _{L_{2}(\Pi)}\rightarrow +\infty$. So there exists a constant $c>0$ such that $\underset{A\in \mathbb{R}^{m_{1}\times m_{2}}}{\min} \psi(A)=\underset{\parallel A\parallel _{L_{2}(\Pi)}\leq c}{\min} \psi(A)$. As the mapping $A\mapsto\rank (A)$ is lower semi-continuous, the functional $\psi (A)$ is lower semi-continuous; thus   $\psi$ attains a minimum on the compact set $\{A:\parallel A\parallel _{L_{2}(\Pi)}\leq c\}$ and the minimum is a global minimum of $\psi$ on $\mathbb{R}^{m_{1}\times m_{2}}$. 

Suppose that Assumption \ref{ass1} is satisfied with equality, i.e., 
$$\parallel A \parallel _{L_{2}(\Pi)}^{2}=\mu^{-2} \parallel A\parallel _{2}^{2}.$$ Then our estimator has a particularly simple form:
\begin{equation}\label{estimator'}
\hat{A}=\underset{A\in \mathbb{R}^{m_{1}\times m_{2}}}{\argmin}\big\{\parallel A-\mathbf{X} \parallel _{2}^{2}+\lambda \mu^{2} \rank (A)\big\},
\end{equation}
where
\begin{equation}
\mathbf{X}=\dfrac{\mu^{2}}{n}\Sum Y_iX_i.
\end{equation}
The optimization problem \eqref {estimator} may equivalently be written as
\begin{equation*}
\hat{A}=\underset{k}{\argmin}\left [\underset{A\in \mathbb{R}^{m_{1}\times m_{2}},\; \rank (A)=k}{\argmin}\parallel A-\mathbf{X} \parallel _{2}^{2}+\lambda \mu^{2} k\right ].
\end{equation*}
Here, the inner minimization problem is to compute the restricted rank estimators $\hat A_k$ that minimizes the norm $\parallel A-\mathbf{X} \parallel _{2}^{2}$ over all matrices of rank $k$. Write the singular value decomposition (SVD) of $\mathbf{X}$:
 \begin{equation}
\mathbf{X}=\overset{\rank\,\mathbf{X}}{\underset{j=1}{\Sigma}}\sigma_j(\mathbf{X})u_j(\mathbf{X})v_j(\mathbf{X})^{T},
\end{equation}
where 
\begin{itemize}
\item  $\sigma_j(\mathbf{X})$ are the singular values of $\mathbf{X}$ indexed in the decreasing order,
\item  $u_j(\mathbf{X})$ (resp. $v_j(\mathbf{X})$) are the left (resp. right) singular vectors of $\mathbf{X}$.
\end{itemize}
Following \cite{reinsel}, one can write:
 \begin{equation}
\hat A_k=\overset{k}{\underset{j=1}{\Sigma}}\sigma_j(\mathbf{X})u_j(\mathbf{X})v_j(\mathbf{X})^{T}.
\end{equation}
Using this, we easily see that $\hat A$ has the form
\begin{equation}
\hat A=\underset{j:\sigma_j(\mathbf{X})\geq\sqrt{\lambda}\mu}{\Sigma}\sigma_j(\mathbf{X})u_j(\mathbf{X})v_j(\mathbf{X})^{T}.
\end{equation}
Thus, the computation of $\hat A$ reduces to hard thresholding of singular values in the SVD of $ \mathbf{X}$.

\textbf{Remark.} We can generalize the estimator given by \eqref{estimator}, taking the minimum over a closed set of the matrices, instead of  the set $\{A\in \mathbb{R}^{m_{1}\times m_{2}}\}$, such as a set of all diagonal matrices, for example.
\section{General oracle inequalities}\label{general}
% % % % % % % % % % % % % % % % % % % % % % % % % % % % % % % % % % % % % % % % % % % % % % % % % % % % %
In the following theorem we bound the rank of our estimator in a very general setting. To the best of our knowledge, such estimates were not known. We also prove general oracle inequalities for the prediction errors in probability analogous to those obtained in \cite[Theorem 2]{Koltchinskii-Tsybakov} for the nuclear norm penalization.

Given $n$ observations $Y_{i}\in\mathbb{R}$ and $X_i$, we define the random matrix
$$ M=\dfrac{1}{n}\overset{n}{\underset{i=1}{\sum}}(Y_{i}X_i-\bE(Y_{i}X_i)).$$
The value $\parallel M\parallel_{\infty}$ determines the ''the noise level'' of our problem.\\ Let $\Delta=\parallel M\parallel_{\infty}$.

\begin{Theorem}\label{thm1}
Let  Assumption \ref{ass1} be satisfied and $\varrho\geq1$. If $\sqrt{\lambda} \geq 2\varrho \mu\Delta$, then
\begin{enumerate}
\item[(i)]\label{rank} \begin{equation*}
\rank(\hat A)\leq \left(1+\dfrac{2}{4\varrho^{2}-1}\right)\rank (A_0),
\end{equation*}
\item  [(ii)]\label{lbound}
\begin{equation*}
\begin{split}
\parallel \hat{A}-A_{0}\parallel _{L_{2}(\Pi)}&\leq \underset{A\in \mathbb{R}^{m_{1}\times m_{2}}}{\mathrm{inf}}\Big \{ \parallel A-A_{0}\parallel _{L_{2}(\Pi)}\\&
 \hskip 0.5 cm+2\sqrt{\lambda\max\left (\dfrac{1}{\varrho^{2}}\rank (A)_{0},\rank (A)\right )}\Big \},
 \end{split}
\end{equation*}
\item[(iii)]
\begin{equation*}
\begin{split}
\parallel \hat A -A_0\parallel _{L_{2}(\Pi)}^{2}&\leq \underset{A\in \mathbb{R}^{m_{1}\times m_{2}}}{\mathrm{inf}}
\Big\{\left (1+\dfrac{2}{2\varrho^{2}-1}\right )\parallel A-A_0 \parallel_{L_{2}(\Pi)}^{2}\\
&\hskip 0.5 cm+2\lambda \left (1+\dfrac{1}{2\varrho^{2}-1}\right )\rank (A)\Big\}.
\end{split}
\end{equation*}
\end{enumerate}
\end{Theorem}
\begin{proof} It follows from the definition of the estimator $\hat A$ that, for all $A\in \matrix$, one has
\begin{equation*}
\begin{split}
&\parallel \hat{A}\parallel _{L_{2}(\Pi)}^{2}-\left \langle \dfrac{2}{n}\Sum Y_iX_i,\hat{A}\right \rangle+\lambda \rank \hat{A}\leq\\&\parallel A\parallel _{L_{2}(\Pi)}^{2}-\left \langle \dfrac{2}{n}\Sum Y_iX_i,A\right \rangle+\lambda \rank (A).
\end{split}
\end{equation*}
Note that
\begin{equation*}
\dfrac{1}{n}\overset{n}{\underset{i=1}{\sum}}\bE(Y_{i}X_i)=\dfrac{1}{n}\overset{n}{\underset{i=1}{\sum}}\bE(\langle A_0,X_i\rangle X_i)
\end{equation*}
and
\begin{equation*}
\dfrac{1}{n}\overset{n}{\underset{i=1}{\sum}}\langle\bE(Y_{i}X_i),A\rangle=
\langle A_0,A\rangle_{L_{2}(\Pi)}.
\end{equation*}
Therefore we obtain
\begin{equation}\label{estim}
\begin{split}
\parallel \hat A -A_0\parallel _{L_{2}(\Pi)}^{2}&\leq\parallel A-A_0 \parallel _{L_{2}(\Pi)}^{2}+2\langle M,\hat A-A\rangle+\lambda (\rank (A)-\rank (\hat A)).
\end{split}
\end{equation}
%\begin{equation*}
%\parallel \hat A \parallel _{L_{2}(\Pi)}^{2}-\langle\dfrac{2}{n}\Sum Y_iX_i,\hat A\rangle+\lambda \rank\hat  A\leq
%\parallel A \parallel _{L_{2}(\Pi)}^{2}-\langle\dfrac{2}{n}\Sum Y_iX_i,A\rangle+\lambda \rank (A)
%\end{equation*}
%Note that $\dfrac{1}{n}\Sum\langle\bE(Y_iX_i),A\rangle=\langle A_0,A\rangle_{L_{2}(\Pi)}$. Therefore, we have
%\begin{equation*}
%\parallel \hat A \parallel _{L_{2}(\Pi)}^{2}-2\langle A_0,\hat A\rangle_{L_{2}(\Pi)}+\lambda \rank\hat  A\leq
%\parallel A \parallel _{L_{2}(\Pi)}^{2}-2\langle A_0,A\rangle_{L_{2}(\Pi)}+2\langle M,\hat A-A\rangle+\lambda \rank (A).
%\end{equation*}
Due to the trace duality $\langle A,B\rangle\leq \parallel  A\parallel_{p}\parallel B\parallel_{q}$ for $p$ and $q$ such that\\ $1/p+1/q=1$, we have
\begin{equation*}
\begin{split}
\parallel \hat A -A_0\parallel _{L_{2}(\Pi)}^{2}&\leq\parallel A-A_0 \parallel _{L_{2}(\Pi)}^{2}+2\Delta\parallel \hat A- A\parallel_{1}+\lambda (\rank (A)-\rank (\hat A))\\
&\leq\parallel A-A_0 \parallel _{L_{2}(\Pi)}^{2}+2\Delta\parallel \hat A- A\parallel_{2}\sqrt{\rank(\hat A- A)}+\\&+\lambda (\rank (A)-\rank (\hat A)).
\end{split}
\end{equation*}
Under Assumption 1, this yields
\begin{equation}\label{triangle}
\begin{split}
\parallel \hat A -A_0\parallel _{L_{2}(\Pi)}^{2}
&\leq\parallel A-A_0 \parallel _{L_{2}(\Pi)}^{2}+\\&+2\mu\Delta\parallel \hat A- A\parallel_{L_{2}(\Pi)}\sqrt{\rank(\hat A- A)}+\lambda (\rank (A)-\rank (\hat A))\\
&\leq\parallel A-A_0 \parallel_{L_{2}(\Pi)}^{2}+\lambda (\rank (A)-\rank (\hat A))\\
&+2\mu\Delta\left(\parallel \hat A-A_0\parallel_{L_{2}(\Pi)}+\parallel A-A_0\parallel_{L_{2}(\Pi)}\right)\sqrt{\rank(\hat A- A)}
\end{split}
\end{equation}
which implies
\begin{equation}\label{2complet}
\begin{split}
&\Big(\parallel \hat A-A_0 \parallel_{L_{2}(\Pi)}-\mu\Delta\sqrt{\rank(\hat A- A)}\Big)^{2}\leq\\
&\leq\Big(\parallel A-A_0 \parallel_{L_{2}(\Pi)}+\mu\Delta\sqrt{\rank(\hat A- A)}\Big)^{2}+\lambda (\rank (A)-\rank (\hat A)).
\end{split}
\end{equation}
To prove (i), we take $A=A_0$ in \eqref{2complet} and we obtain:
\begin{equation}\label{rank1}
\begin{split}
\lambda(\rank (\hat A) -\rank (A_0))&\leq \Big(\mu\Delta\sqrt{\rank(\hat A- A_0)}\Big)^{2}\\
&\leq \dfrac{\lambda}{4\varrho^2}(\rank (\hat A) +\rank (A_0)).
\end{split}
\end{equation}
Thus,
\begin{equation}
\rank(\hat A)\leq \left(1+\dfrac{2}{4\varrho^{2}-1}\right)\rank (A_0).
\end{equation} 
To prove (ii), we first consider the case $\rank (A)\leq \rank (\hat A)$. Then \eqref{2complet} implies
\begin{equation*}
\begin{split}
0\leq&\lambda (\rank (\hat A)-\rank (A))\leq \Big(\parallel A-A_0 \parallel_{L_{2}(\Pi)}+\parallel \hat A-A_0 \parallel_{L_{2}(\Pi)}\Big)\times\\&\Big(\parallel A-A_0 \parallel_{L_{2}(\Pi)}-\parallel \hat A-A_0 \parallel_{L_{2}(\Pi)}+2\mu\Delta\sqrt{\rank(\hat A- A)}\Big).
\end{split}
\end{equation*}
Therefore, for $\rank\,A\leq\rank(\hat A)$, we have
\begin{equation}\label{cas1}
\begin{split}
\parallel \hat A-A_0 \parallel_{L_{2}(\Pi)}&\leq\parallel A-A_0 \parallel_{L_{2}(\Pi)}+2\mu\Delta\sqrt{\rank(\hat A- A)}\\
&\leq\parallel A-A_0 \parallel_{L_{2}(\Pi)}+2\mu\Delta\sqrt{\rank (\hat A)+\rank (A)}\\
&\leq\parallel A-A_0 \parallel_{L_{2}(\Pi)}+\sqrt{\dfrac{2\lambda}{\varrho^{2}} \rank (\hat A)}.
\end{split}
\end{equation}
Using (i) we obtain
\begin{equation}
\label{lb1}
\begin{split}
\parallel \hat A-A_0 \parallel_{L_{2}(\Pi)}
\leq\parallel A-A_0 \parallel_{L_{2}(\Pi)}+2\sqrt{\dfrac{\lambda}{\varrho^{2}} \rank (A_0)}.
\end{split}
\end{equation}

Consider now the case, $\rank (A) \geq \rank (\hat A)$. Using that $\sqrt{a^{2}+b^2}\leq a+b$ for $a\geq0$ and $b\geq0$, we get  from \eqref{2complet} that
\begin{equation}\label{cas2}
\begin{split}
\parallel \hat A-A_0 \parallel_{L_{2}(\Pi)}&\leq\parallel A-A_0 \parallel_{L_{2}(\Pi)}+2\mu\Delta\sqrt{\rank(\hat A- A)}+\\&+\sqrt{\lambda (\rank (A)-\rank (\hat A))}\\&\leq \parallel A-A_0 \parallel_{L_{2}(\Pi)}+\sqrt{\lambda}\left (\sqrt{\rank(\hat A)+\rank (A)}+\sqrt{\rank (A)-\rank (\hat A)}\right ).
\end{split}
\end{equation}
Finally, the elementary inequality $\sqrt{a+c}+\sqrt{a-c}\leq 2\sqrt{a}$ yields
\begin{equation}\label{lb2}
\parallel \hat A-A_0 \parallel_{L_{2}(\Pi)}\leq \parallel A-A_0 \parallel_{L_{2}(\Pi)}+2\sqrt{\lambda\rank (A)}.
\end{equation}
Using \eqref{lb1} and \eqref{lb2}, we obtain (ii).

To prove (iii), we use \eqref{triangle} to obtain
\begin{eqnarray*}
\parallel \hat A -A_0\parallel _{L_{2}(\Pi)}^{2}&\leq&
\parallel A-A_0 \parallel_{L_{2}(\Pi)}^{2}+\lambda (\rank (A)-\rank (\hat A))+\\
&+&2\dfrac{\sqrt{\frac{\lambda}{2}}}{\varrho\sqrt{2}}\parallel \hat A-A_0\parallel_{L_{2}(\Pi)}\sqrt{\rank(\hat A)+ \rank (A)}\\
&+&2\dfrac{\sqrt{\frac{\lambda}{2}}}{\varrho\sqrt{2}}\parallel A-A_0\parallel_{L_{2}(\Pi)}\sqrt{\rank(\hat A)+\rank (A)}.
\end{eqnarray*}
From which we get
\begin{equation*}
\begin{split}
\left (1-\dfrac{1}{2\varrho^{2}}\right )\parallel \hat A -A_0\parallel _{L_{2}(\Pi)}^{2}&\leq
\left (1+\dfrac{1}{2\varrho^{2}}\right )\parallel A-A_0 \parallel_{L_{2}(\Pi)}^{2}\\&+\lambda( \rank (A) +\rank(\hat A))+\lambda( \rank (A) -\rank(\hat A))
\\&\leq\left (1+\dfrac{1}{2\varrho^{2}}\right )\parallel A-A_0 \parallel_{L_{2}(\Pi)}^{2}+2\lambda\, \rank (A)
\end{split}
\end{equation*}
and (iii) follows.
\end{proof} 
In the next theorem we obtain bounds for the prediction error in expectation. Set $m=m_1+m_2$, $m_1\wedge m_2=\min (m_1,m_2)$ and $m_1\vee m_2=\max (m_1,m_2)$. Suppose that $\bE\left (\Delta^{2}\right )<\infty$ and let $\mathcal B_r$ be the set of non-negative random variables $W$ bounded by $r$. We set 
\begin{equation*} %\label{estEnormMrank-1}
 \begin{split}
 S=\underset{W\in\mathcal B_{m_1\wedge m_2}}{\sup}\dfrac{\bE\left (\Delta^{2}W\right)}{ \max\{\bE(W), 1\}} \leq\left ( m_1\wedge m_2\right )\bE\left (\Delta^{2}\right )<\infty.
 \end{split}
 \end{equation*}
% % % % % % % % % % % % % % % % % % % % % % % % % % % % % % % % % % % % % % % % % % % % % % % % % % % % % % % % % % % %
\begin{Theorem}\label{thm2}
Let Assumption \ref{ass1} be satisfied. Consider $\varrho\geq 1$ and a regularization parameter $\lambda$ satisfying $\sqrt{\lambda}\geq 2\varrho\mu\,\sqrt{S}$. Then
\begin{enumerate}
\item [(a)] 
\begin{equation*}
\begin{split}
\bE(\rank(\hat A))\leq \max\left \{ \left (1+\dfrac{2}{4\varrho^{2}-1}\right )\rank (A_0), \dfrac{1}{4\varrho^{2}}\right \},
\end{split}
\end{equation*}
\item [(b)] 
\begin{equation*}
\begin{split}
\bE\left (\parallel \hat A-A_0 \parallel_{L_{2}(\Pi)}\right )&\leq\underset{A\in \mathbb{R}^{m_{1}\times m_{2}}}{\mathrm{inf}}\Big\{\parallel A-A_0 \parallel_{L_{2}(\Pi)}\\&\hskip 0.5 cm+\dfrac{5}{2}\sqrt{\lambda\max\left(\rank (A), \dfrac{\rank (A_0)}{\varrho^{2}},\dfrac{1}{4\varrho^{2}} \right)}\Big\},
\end{split}
\end{equation*}
and
\item [(c)] 
\begin{equation*}
\begin{split}
\bE\left (\parallel \hat A-A_0 \parallel_{L_{2}(\Pi)}^{2}\right )&\leq\underset{A\in \mathbb{R}^{m_{1}\times m_{2}}}{\mathrm{inf}}\left \{ \left (1+\dfrac{2}{2\varrho^{2}-1}\right )\parallel A-A_0 \parallel_{L_{2}(\Pi)}^{2}\right .\\&\hskip 0.5 cm\left .+2\lambda\left (1+\dfrac{1}{2\varrho^{2}-1}\right )\max\left (\rank (A),\dfrac{1}{2}\right )\right \}.
\end{split}
\end{equation*}
\end{enumerate}
\end{Theorem}
\begin{proof}
To prove (a)  we  take the expectation of \eqref{rank1} to obtain
\begin{equation}\label{E(a)}
\begin{split}
\lambda(\bE(\rank (\hat A)) -\rank (A_0))\leq \bE \Big(\mu^{2}\Delta^{2}\left (\rank(\hat A) +\rank\,A_0\right )\Big).
\end{split}
\end{equation}
If $A_0=0$, as $\sqrt{\lambda}\geq 2\varrho\mu\,C$ we obtain
\begin{equation}
\lambda\bE(\rank (\hat A))\leq  \mu^{2}C^{2}  \max\{\bE(\rank (\hat A)), 1\}\leq \dfrac{\lambda}{4\varrho^{2}} \max\{\bE(\rank (\hat A)), 1\}
\end{equation}
which implies $\bE(\rank (\hat A))\leq \dfrac{1}{4\varrho^{2}}$.

If $A_0\neq 0$, $\rank(\hat A)+\rank (A_0)\geq 1$ and we get
\begin{equation*}
\begin{split}
\lambda(\bE(\rank (\hat A)) -\rank (A_0))&\leq \mu^{2}\,C^{2}\left (\bE \left (\rank(\hat A)\right)+\rank (A_0)\right)\\&\leq \dfrac{\lambda}{4\,\varrho^{2}}\left (\bE \left (\rank(\hat A)\right)+\rank (A_0)\right)
\end{split}
\end{equation*}
which proves part (a) of Theorem \ref{thm2}.
%\begin{equation}\label{estEranklambda}
%\begin{split}
%\lambda(\bE(\rank (\hat A)) -\rank (A_0))\leq \dfrac{\lambda}{4\,\varrho^{2}}\left (\left (\bE \left (\rank(\hat A)\right )^{q}\right )^{1/q}+\rank (A_0)\right ).
%\end{split}
%\end{equation}
%As $\rank\; \hat A\leq m_1\wedge m_2$ we compute
%\begin{equation*}
%\left (\rank (\hat A)\right )^{1/(log (m)-1)}\leq \exp (\log(m_1\wedge m_2)/(log (m)-1))\leq e^{2}
%\end{equation*}
%and as $\rank (\hat A)\geq 1$ we obtain
%\begin{equation}\label{estErank^q}
%\left (\bE\left (\rank (\hat A)\right )^{1+1/(log (m)-1)}\right )^{1-1/\log m}\leq \bE\left (\left (\rank (\hat A)\right )^{1+1/(log (m)-1)}\right )\leq e^{2}\;\bE( \rank (\hat A)).
%\end{equation}
%Plugging \eqref{estErank^q} into \eqref{estEranklambda} we get (a).

To prove (b), \eqref{cas1} and \eqref{cas2} yield
\begin{equation*} %\label{bestim1}
\begin{split}
\parallel \hat A-A_0 \parallel_{L_{2}(\Pi)}&\leq\parallel A-A_0 \parallel_{L_{2}(\Pi)}+2\mu\Delta\sqrt{\rank (\hat A)+\rank (A)}\\&\hskip 0.5 cm+\mathbb{I}_{\rank (\hat A)\leq \rank (A)}\sqrt{\lambda(\rank (A)-\rank (\hat A))}
\end{split}
\end{equation*}
where $\mathbb{I}_{\rank (\hat A)\leq \rank (A)}$ is the indicator function of the event $\{\rank (\hat A)\leq \rank (A)\}$.
Taking the expectation we obtain
\begin{equation*} %\label{bestim1}
\begin{split}
\bE\parallel \hat A-A_0 \parallel_{L_{2}(\Pi)}&\leq\parallel A-A_0 \parallel_{L_{2}(\Pi)}+2\mu\bE\left (\Delta\sqrt{\rank (\hat A)+\rank (A)}\right )\\&\hskip 0.5 cm+\sqrt{\lambda}\bE\left (\mathbb{I}_{\rank (\hat A)\leq \rank (A)}\sqrt{\rank (A)-\rank (\hat A)}\right) .
\end{split}
\end{equation*}
Note that  Cauchy-Schwarz inequality and $C^{2}\geq S$ imply
\begin{equation*} % \label{estEnormMrank-2}
 \begin{split}
 \bE\left (\Delta W\right) \leq C\max\{\bE(W), 1\}.
 \end{split}
 \end{equation*}
Taking $W=\sqrt{\rank (\hat A)+\rank (A)}$ we find
\begin{equation} \label{bestim2}
\begin{split}
\bE\parallel \hat A-A_0 \parallel_{L_{2}(\Pi)}&\leq\parallel A-A_0 \parallel_{L_{2}(\Pi)}+2\mu\,C \max\left \{1,\bE\sqrt{\rank (\hat A)+\rank (A)}\right \}\\&\hskip 0.5 cm+\sqrt{\lambda}\bE\left (\mathbb{I}_{\rank (\hat A)\leq \rank (A)}\sqrt{\rank (A)-\rank (\hat A)}\;\right) .
\end{split}
\end{equation}
If $\bE\sqrt{\rank (\hat A)+\rank (A)}<1$, which implies $A=0$, as $\sqrt{\lambda}\geq 2\varrho\mu\,C$ we obtain
\begin{equation*} %\label{bestim2}
\begin{split}
\bE\parallel \hat A-A_0 \parallel_{L_{2}(\Pi)}&\leq\parallel A-A_0 \parallel_{L_{2}(\Pi)}+\dfrac{\sqrt{\lambda}}{\varrho}.
\end{split}
\end{equation*}
This prove (b) in the case $\bE\sqrt{\rank (\hat A)+\rank (A)}<1$. 

If $\bE\sqrt{\rank (\hat A)+\rank (A)}\geq 1$, from \eqref{bestim2} we get
\begin{equation*} %\label{bestim2}
\begin{split}
\bE\parallel \hat A-A_0 \parallel_{L_{2}(\Pi)}\leq\parallel A-A_0 \parallel_{L_{2}(\Pi)}&+\sqrt{\lambda}\left \{\dfrac{1}{\varrho}\bE\left (\mathbb{I}_{\rank (\hat A)\leq \rank (A)}\sqrt{\rank (\hat A)+\rank (A)}\right )\right .\\&\hskip 1 cm\left .+\bE\left (\mathbb{I}_{\rank (\hat A)\leq \rank (A)}\sqrt{\rank (A)-\rank (\hat A)}\right)\right .\\&\hskip 1 cm\left .+ \dfrac{1}{\varrho}\bE\left (\mathbb{I}_{\rank (\hat A)> \rank (A)}\sqrt{\rank (\hat A)+\rank (A)}\right )\right\}.
\end{split}
\end{equation*}
Using that $\varrho\geq 1$ and the elementary inequality $\sqrt{a+c}+\sqrt{a-c}\leq 2\sqrt{a}$  we find
\begin{equation*} %\label{bestim2}
\begin{split}
\bE\parallel \hat A-A_0 \parallel_{L_{2}(\Pi)}&\leq\parallel A-A_0 \parallel_{L_{2}(\Pi)}+\sqrt{\lambda}\Big (2\sqrt{\rank (A)}\;\mathbb{P}(\rank (\hat A)\leq \rank (A))\\&\hskip 0.5 cm+ \dfrac{1}{\varrho}\bE\big(\mathbb{I}_{\rank (\hat A)> \rank (A)}\sqrt{2\rank (\hat A)}\big)\Big ).
\end{split}
\end{equation*}
The Cauchy-Schwarz inequality and (a)  imply
\begin{equation*} %\label{bestim2}
\begin{split}
\bE\parallel \hat A-A_0 \parallel_{L_{2}(\Pi)}&\leq\parallel A-A_0 \parallel_{L_{2}(\Pi)}+\sqrt{\lambda}\Big (2\sqrt{\rank (A)}\;\mathbb{P}(\rank (\hat A)\leq \rank (A))\\&+ \dfrac{2}{\varrho}\sqrt{\max\left \{\rank (A_0),\dfrac{1}{4\varrho^{2}}\right \}}\;\mathbb{P}^{1/2}(\rank (\hat A)> \rank (A))\Big ).
\end{split}
\end{equation*}
Using that $x+\sqrt{1-x}\leq 5/4$ when $0\leq x\leq 1$ for $x=\mathbb{P}(\rank (\hat A)\leq \rank (A))$ we get (b).

We now prove part (c). From \eqref{triangle} we compute
\begin{equation*} %\label{triangle}
\begin{split}
\parallel \hat A -A_0\parallel _{L_{2}(\Pi)}^{2}
&\leq\parallel A-A_0 \parallel_{L_{2}(\Pi)}^{2}+\lambda (\rank (A)-\rank (\hat A))\\&
+2\left (\sqrt{2}\mu\varrho\Delta\sqrt{\rank(\hat A)+ \rank (A)}\right )\times\\
&\hskip 0.5 cm\times\left(\dfrac{1}{\sqrt{2}\varrho}\parallel \hat A-A_0\parallel_{L_{2}(\Pi)}+\dfrac{1}{\sqrt{2}\varrho}\parallel A-A_0\parallel_{L_{2}(\Pi)}\right)\\
&\leq \left (1+\dfrac{1}{2\varrho^{2}}\right ) \parallel A-A_0 \parallel_{L_{2}(\Pi)}^{2}+\dfrac{1}{2\varrho^{2}}\parallel \hat A-A_0 \parallel_{L_{2}(\Pi)}^{2}\\
&+4\varrho^{2}\mu^{2}\Delta^{2}(\rank(\hat A)+\rank (A))+\lambda (\rank (A)-\rank (\hat A))\\
\end{split}
\end{equation*}
which implies
\begin{equation*}
\begin{split}
\left (1-\dfrac{1}{2\varrho^{2}}\right )&\parallel \hat A-A_0 \parallel_{L_{2}(\Pi)}^{2}\leq\left (1+\dfrac{1}{2\varrho^{2}}\right )\parallel A-A_0 \parallel_{L_{2}(\Pi)}^{2}\\ &+4\;\varrho^{2}\mu^{2}\Delta^{2}(\rank(\hat A)+\rank (A))+\lambda(\rank (A)-\rank (\hat A)) .
\end{split}
\end{equation*}
Taking the expectation we obtain 
\begin{equation*} %\label{th2iiE}
\begin{split}
\bE\parallel \hat A-A_0 \parallel_{L_{2}(\Pi)}^{2}&\leq\left (1-\dfrac{1}{2\varrho^{2}}\right )^{-1}\Big\{\left (1+\dfrac{1}{2\varrho^{2}}\right )\parallel A-A_0 \parallel_{L_{2}(\Pi)}^{2}\\&\hskip 0.5 cm+4\;\varrho^{2}\mu^{2}\bE\left (\Delta^{2}(\rank(\hat A)+\rank (A))\right )\\&\hskip 0.5 cm+\lambda\bE(\rank (A)-\rank (\hat A))\Big\} .
\end{split}
\end{equation*}
%By H{\"o}lder's inequality and \eqref{estEM} we get
%\begin{equation*}
%\begin{split}
%\bE\left (\Delta^{2}\rank(\hat A)\right)&\leq \left (\bE \Delta^{2\log m}\right )^{1/\log m} \left (\bE(\rank(\hat A))^{q}\right )^{1/q}\\&\hskip 0.5 cm\leq \dfrac{\log m}{\nu^{2}} \left (\bE(\rank(\hat A))^{q}\right )^{1/q}.
%\end{split}
%\end{equation*}
As $C^{2}\geq S$ we compute
\begin{equation}\label{final}
\begin{split}
\bE\parallel \hat A-A_0 \parallel_{L_{2}(\Pi)}^{2}&\leq\left (1-\dfrac{1}{2\varrho^{2}}\right )^{-1}\Big\{\left (1+\dfrac{1}{2\varrho^{2}}\right )\parallel A-A_0 \parallel_{L_{2}(\Pi)}^{2}\\&+4\varrho^{2}\mu^{2}\,C^{2}\max \left(1,\bE(\rank(\hat A))+\rank (A)\right  )\\&+\lambda(\rank (A)-\bE(\rank (\hat A)))\Big\} .
\end{split}
\end{equation}
The assumption on $\lambda$ and \eqref{final} imply (c).
%\begin{equation*}
%\begin{split}
%\bE\parallel \hat A-A_0 \parallel_{L_{2}(\Pi)}^{2}&\leq\left (1-\dfrac{e^{2}}{2\varrho^{2}}\right )^{-1}\left \{\left (1+\dfrac{e^{2}}{2\varrho^{2}}\right )\parallel A-A_0 \parallel_{L_{2}(\Pi)}^{2}\right .\\&\left .+\lambda\Big (  \dfrac{1}{e^{2}}\big(\bE(\rank(\hat A))^{q}\big )^{1/q}-\bE(\rank (\hat A))+(1+1/e^{2})\rank (A)\Big)\right \}\\&\leq \left (1+\dfrac{e^{2}}{2\varrho^{2}}\right )\parallel A-A_0 \parallel_{L_{2}(\Pi)}^{2}+\lambda\left (1+1/e^{2}\right )\rank (A).
%\end{split}
%\end{equation*}
This completes the proof of  Theorem \ref{thm2}.
\end{proof}
 The next lemma gives an upper bound on $S$ in the case when $ \Delta$ concentrates exponentially around its mean. 
  \begin{lemma}\label{concentratio_bound}
           Assume that
           \begin{equation}\label{distr_lemma_conc}
            \mathbb P\left \{\Delta\geq \bE \Delta + t\right \}\leq \exp\left \{-ct^{\alpha}\right \}.
            \end{equation}
            for some positive constants $c$ and $\alpha$. Then
            \begin{equation}\label{est_S}
           S\leq 2\left (\bE \Delta\right)^{2}+e\,2^{1+1/p}\left (p/c\alpha\right )^{2/\alpha}
           \end{equation}
           for $p\geq \max\{2\log(m_1\wedge m_2)+1,\alpha\}$.
           \end{lemma}
           \begin{proof}
           Write
           \begin{equation*}
      \begin{split}
       \bE\left (\Delta^{2}W\right)&=\left (\bE\left (\Delta\right )\right )^{2}\bE (W)+
      \bE\left [\Big (\Delta-\bE \left (\Delta\right )\Big )^{2}W\right]
      \\&+2\bE\left (\Delta\right )\bE\left [\Big (\Delta-\bE \left (\Delta\Big )\right )W\right]\\
      &\leq 2\left (\bE\left (\Delta\right )\right )^{2}\max\{\bE( W),1\}\\&+ \bE\left [\Big (\Delta-\bE \left (\Delta\right )\Big )^{2}W\right]\\&+
       \bE\left [\Big (\Delta-\bE \left (\Delta\right )\Big )^{2}W^{2}\right].
      \end{split}
           \end{equation*}
           Setting $\bar W =\max\{W,W^{2}\}$, we see that it is enough to estimate $\bE\left [\Big (\Delta-\bE \left (\Delta\right )\Big )^{2}\bar W\right]$ for $0\leq \bar W\leq (m_1\wedge m_2)^{2}$. Putting $X=\Delta-\,\bE \left (\Delta\right )$
            and using H{\"o}lder's inequality we get
           \begin{equation}\label{estEnormMrank}
           \begin{split}
           \bE \left (X^{2}\bar W\right)\leq \left (\bE X^{2p}\right )^{1/p} \left (\bE\,\bar W^{q}\right )^{1/q}
           \end{split}
           \end{equation}
           where $q=1+1/(p-1)$.
           We first estimate $\left (\bE X^{2p}\right )^{1/p}$. Inequality \eqref{distr_lemma_conc} implies that 
           \begin{equation}\label{estEM-lc}
           \begin{split}
           \left (\bE X^{2p}\right )^{1/p} &=\left ( \overset{+\infty}{\underset{0}{\int}}\mathbb P \left (X> t^{1/(2p)}\right )\mathrm{d}t\right )^{1/p}
          \leq
         \left ( \overset{+\infty}{\underset{0}{\int}}\exp\{-ct^{\alpha/(2p)}\}\mathrm{d}t
      \right )^{1/p}
\\&= c^{-2/\alpha}\left (\dfrac{2p}{\alpha}\Gamma\left (\dfrac{2p}{\alpha}\right )\right )^{1/p}.
           %\\&=\left ( \overset{+\infty}{\underset{0}{\int}} 2\log(m)\; m\; t^{2\log( m)-1} e^{-t^{2}\nu^{2}} \mathrm{d}t\right )^{1/\log m}
           \end{split}
           \end{equation}
           The Gamma-function satisfies the following bound: 
           \begin{equation}\label{Gamma}
           \text{for}\quad x\geq 2,\quad\Gamma(x)\leq \left (\dfrac{x}{2}\right )^{x-1},
           \end{equation} 
           cf. Proposition \ref{proposition1}. Plugging it into \eqref{estEM-lc} we find  
           \begin{equation}\label{X_p}
           \begin{split}
           \left (\bE X^{2p}\right )^{1/p}& \leq 2^{1/p}\left (\dfrac{p}{c\alpha}\right )^{2/\alpha}.
           \end{split}
           \end{equation}
           If $\bE(\bar W^{q})<1$ we get \eqref{est_S} directly from \eqref{X_p}.
           If $\bE(\bar W^{q})\geq 1$, the bound $\bar W\leq \left (m_1\wedge m_2\right )^{2}$ implies that
           \begin{equation*} %\label{W_q}
           \bar W^{1/(p-1)}\leq e
           \end{equation*}
           and thus
           \begin{equation}\label{W_q}
           \left (\bE \left (\bar W^{1+1/(p-1)}\right )\right )^{1-1/p}\leq e\bE(\bar W).\end{equation}
           Then \eqref{est_S} follows from \eqref{X_p} and \eqref{W_q}.
                     \end{proof}
 
% % % % % % % % % % % % % % % % % % % % % % % % % % % % % % % % % % % % % % % % % % % % % % % % % % % % % % %
%%%%%%%%%%%%%%%%%%%%%%%%%%%%%%%%%%%%%%%%%%%%%%%%%%%%%%%%%%%%%%%%%%%%%%%%%%%%%%%%%%%%%%%%%%%%%%%%%%%%%%%%%%%%%%%%%
\section{Matrix Completion}\label{matrix}
In this section we present some consequences of the general oracles inequalities of Theorems \ref{thm1} and \ref{thm2} for the model of USR matrix completion. Assume that the design matrices $X_i$ are i.i.d uniformly distributed on the set $\mathcal{X}$ defined in \eqref{basisUSR}. This implies that
\begin{equation}
m_1m_2\parallel A\parallel_{L_{2}(\Pi)}^{2}=\parallel A \parallel^{2}_{2},
\end{equation}
for all matrices $A\in\matrix$.
%\begin{equation}\label{estimator'}
%\hat{A}=\underset{A\in \mathbb{R}^{m_{1}\times m_{2}}}{\argmin}\big\{\parallel A-\mathbf{X} \parallel _{2}^{2}+\lambda m_1m_2 \rank (A)\big\},
%\end{equation}
%where
%\begin{equation}
%\mathbf{X}=\dfrac{m_1m_2}{n}\Sum Y_iX_i.
%\end{equation}
%The optimization problem \eqref {estimator'} may equivalently be written as
%\begin{equation*}
%\hat{A}=\underset{k}{\argmin}\underset{A\in \mathbb{R}^{m_{1}\times m_{2}},\; \rank (A)=k}{\argmin}\big[\parallel A-\mathbf{X} \parallel _{2}^{2}+\lambda m_1m_2 k\big].
%\end{equation*}
%One needs to compute the restricted rank estimators $\hat A_k$ that minimize the norm $\parallel A-\mathbf{X} \parallel _{2}^{2}$ over all matrices of rank $k$. Write the singular value decomposition (SVD) of $X$
% \begin{equation}
%\hat A_k=\overset{\rank\,X}{\underset{j=1}{\Sigma}}\sigma_j(\mathbf{X})u_j(\mathbf{X})v_j(\mathbf{X})^{T}
%\end{equation}
%where 
%\begin{itemize}
%\item  $\sigma_j(\mathbf{X})$ are the singular values of $\mathbf{X}$ indexed in the decreasing order,
%\item  $u_j(\mathbf{X})$ (resp. $v_j(\mathbf{X})$) are the left (resp. right) singular vectors of $\mathbf{X}$.
%\end{itemize}
%Following \cite{reinsel}, one can write:
% \begin{equation}
%\hat A_k=\overset{k}{\underset{j=1}{\Sigma}}\sigma_j(\mathbf{X})u_j(\mathbf{X})v_j(\mathbf{X})^{T}.
%\end{equation}
Then, we can write $\hat A$ explicitly
\begin{equation}\label{hardtresholding}
\hat A=\underset{j:\sigma_j(\mathbf{X})\geq\sqrt{\lambda m_1m_2}}{\Sigma}\sigma_j(\mathbf{X})u_j(\mathbf{X})v_j(\mathbf{X})^{T}.
\end{equation}
Set $\hat{r}=\rank(\hat A)$. In the case of matrix completion, we can improve point (i) of Theorem \ref{thm1} and give an estimation on the difference of the first $\hat r$ singular values of $\hat A$ and $A_0$. We also get  bounds on the prediction error measured in norms different from the Frobenius norm, in particular in the spectral norm.
\begin{Theorem}\label{thm3}
Let $\lambda$ satisfy the inequality $\sqrt{\lambda}\geq 2\mu \Delta$ (as in Theorem \ref{thm1}). Then 
\begin{enumerate}
\item [(i)] $\hat r\leq \rank(A_0)$;
\item [(ii)] $\mid\sigma_j(\hat A)-\sigma_j(A_0)\mid\leq \dfrac{\sqrt{\lambda m_1m_2}}{2}$ for $j=1,\dots,\hat r$;
\item[(iii)] $ \left \|\hat A-A_0\right \|_{\infty}\leq \dfrac{3}{2}\sqrt{\lambda m_1m_2}$;
\item[(iv)] for $2\leq q\leq \infty$, one has $$\left \|\hat A-A_0\right \|_q\leq \dfrac{3}{2}\left (4/3\right )^{2/q}\sqrt{m_1m_2\lambda}(\rank (A_0))^{1/q},$$
where we set $x^{1/q}=1$ for $x>0,q=\infty$.
\end{enumerate}
\end{Theorem} 
\begin{proof} The proof is obtained by adapting the proof of \cite[Theorem 8]{Koltchinskii-Tsybakov} to hard thresholding estimators. For completeness, we give the proof of (iii) and (iv). 

Let us start with the proof of (iii). Note that $\mathbf{X}-A_0=m_1m_2M$. Let $B=\mathbf{X}-\hat A$, by \eqref{hardtresholding} we have that $\sigma_1(B)\leq \sqrt{\lambda m_1m_2}$. Then
\begin{equation}
\begin{split}
\sigma_1(\hat A-A_0)&=\sigma_1(\mathbf{X}-A_0-B)=\sigma_1(m_1m_2M-B)\\&\leq m_1m_2 \Delta+\sqrt{\lambda m_1m_2}
\leq \dfrac{3}{2}\sqrt{\lambda m_1m_2}
\end{split}\end{equation}
and we get (iii).

To prove (iv) we use the following interpolation inequality (see \cite[Lemma 11]{2009arXiv0912.5338R}):
for $0<p<q<r\leq \infty$ let $\theta\in [0,1]$ be such that $\dfrac{\theta}{p}+\dfrac{1-\theta}{r}=\dfrac{1}{q}$ then for all $A\in \mathbb{R}^{m_{1}\times m_{2}}$we have

\begin{equation}\label{interpolation}
\left \|A\right \|_q\leq  \left \|A\right \|^{\theta}_p \left \|A\right \|^{1-\theta}_r.
\end{equation}
For $q\in (2,\infty)$ take $p=2$ and $r=\infty$. From Theorem \ref{thm1} (ii) we get that
\begin{equation}\label{normS2}
 \left \|\hat A-A_0\right \|_2\leq 2\sqrt{m_1m_2\lambda\rank (A_0)}.
\end{equation}
Now, plugging (iii) of Theorem \ref{thm3} and \eqref{normS2} into \eqref{interpolation}, we obtain
\begin{equation}
 \left \|\hat A-A_0\right \|_q\leq \left \|\hat A-A_0\right \|^{2/q}_2\left \|\hat A-A_0\right \|_\infty^{1-2/q}\leq 
 \dfrac{3}{2}\left (4/3\right )^{2/q}\sqrt{m_1m_2\lambda}(\rank (A_0))^{1/q}
\end{equation}
and (iv) follows. This completes the proof of Theorem \ref{thm3}.
\end{proof}

 In view of Theorems \ref{thm1} and \ref{thm2}, to specify the value of regularization parameter $\lambda$, we need to estimate $\Delta$ with high probability. We will use the bounds obtained in \cite{Koltchinskii-Tsybakov} in the following two settings of particular interest:
 \begin{itemize}
\item [(\textbf{A})] \textit{Statistical learning setting}. 
There exists a constant $\eta$ such that\\ $\underset{i=1,\dots, ,n}{\max}\mid Y_i\mid\leq \eta$. Then, we set 
\begin{equation}\label{ro1}
\rho(m_1,m_2,n,t)=4\eta\max \left\{\sqrt{\dfrac{t+\log(m)}{(m_1\wedge m_2)n}},\dfrac{2(t+\log(m))}{n}\right \},
\end{equation}
\begin{equation*}
n^{*}=4(m_1\wedge m_2)\log m,\qquad c^{*}=4\eta.
\end{equation*}
\item [(\textbf B)] \textit{Sub-exponential noise.} We suppose that the pairs $(X_i,Y_i)_{i}$ are iid and that there exist constants $\omega,c_1>0,\alpha\geq 1$ and $c_2$ such that
\begin{equation*}
\underset{i=1,\dots, ,n}{\max}\bE\exp\left (\dfrac{\mid\xi_i\mid^{\alpha}}{\omega^{\alpha}}\right )<c_2,\quad\bE\xi_i^{2}\geq c_1\omega^{2},\;\forall\:1\leq i\leq n.
\end{equation*}
Let $A_0=(a^{0}_{ij})$ and $\underset{i,j}{\max}\mid a^{0}_{ij}\mid \leq a$. Then, we set 
\begin{equation}\label{ro2}
\begin{split}
\rho(m_1,m_2,n,t)&=\tilde C(\omega\vee a)\max \left\{\sqrt{\dfrac{t+\log(m)}{(m_1\wedge m_2)n}},\right .\\&\hskip 0.5 cm\left .\dfrac{(t+\log(m))\log^{1/\alpha}(m_1\wedge m_2)}{n}\right \},
\end{split}
\end{equation}
\begin{equation*}
n^{*}=(m_1\wedge m_2)\log^{1+2/\alpha} (m),\qquad c^{*}=\tilde C (\omega\vee a). 
\end{equation*}
where $\tilde C>0$ is a large enough constant that depends only on $\alpha,c_1,c_2$.
\end{itemize}
In both case we can estimate  $\Delta$ with high probability:
\begin{lemma}[\cite{Koltchinskii-Tsybakov}, Lemmas 1, 2 and 3]\label{lemma1} 
For all $t>0$, with probability at least $1-e^{-t}$ in the case of statistical learning setting (respectively, $1-3e^{-t}$ in the case of sub-exponential noise), one has
\begin{equation}\label{estTsyblemma}
\Delta\leq \rho(m_1,m_2,n,t). 
\end{equation}
\end{lemma}
As a corollary of Lemma \ref{lemma1} we obtain the following bound for $$ S=\underset{\mathcal B_{m_1\wedge m_2}}{\sup}\dfrac{\bE\left (\Delta^{2}W\right)}{ \max\{\bE(W), 1\}}.$$
\begin{lemma}\label{lemma2}
Let one of the set of conditions (\textbf A) or (\textbf B) be satisfied. Assume $n>n^{*}$, $\log m\geq 5$ and $W$ is a non-negative random variable such that $W\leq m_1\wedge m_2$,  then
 \begin{equation}\label{estEnormMrank-1}
 \begin{split}
 \bE\left (\Delta^{2}W\right) \leq \dfrac{(c^{*}e)^{2}\log m}{n(m_1\wedge m_2)}  \max\{\bE(W), 1\}.
 \end{split}
 \end{equation}
% and
% \begin{equation}\label{estEnormMrank-2}
%  \begin{split}
%  \bE\left (\DeltaW\right)\leq \dfrac{e\sqrt{\log m}}{\nu} \max\{\bE(W), 1/\sqrt{e}\}.
%  \end{split}
%  \end{equation}}
\end{lemma}
\begin{proof}We will prove \eqref{estEnormMrank-1} in the case of statistical learning setting. The proof in the case of sub-exponential noise is completely analogous.  Set $$t^{*}=\dfrac{n}{4(m_1\wedge m_2)}-\log m.$$ % for the statistical learning setting and $$t^{*}=\dfrac{n}{(m_1\wedge m_2)\log^{2/\alpha}((m_1\wedge m_2)}-\log m$$
%in the case of sub-exponential noise.
Note that Lemma \ref{lemma1} implies that 
\begin{equation}\label{proba_1}
\mathbb P \left (\Delta> t\right )\leq m\exp\{-t^{2}\,n(m_1\wedge m_2)(c^{*})^{-2}\}\qquad\text{for}\qquad t\leq t^{*}
\end{equation}
 and
\begin{equation}\label{proba_2}
\mathbb P \left (\Delta> t\right )\leq3\sqrt{m}\exp\{-t\,n/(2c^{*})\}\qquad\text{for}\qquad t\geq t^{*}.
\end{equation}
 % and statistical learning setting, and
%\begin{equation}\label{proba_}
%\mathbb P \left (\Delta> t\right )\leq m \exp\{-t\,n/(c^{*}\log^{1/\alpha}(m_1\wedge m_2))\}.
%\end{equation}
%for $t\geq t^{*}$ and sub-exponential noise.}

We set $\nu_1=n(m_1\wedge m_2)(c^{*})^{-2}$, $\nu_2=n(2c^{*})^{-1}$ %, $\nu_3=n/(c^{*}\log^{1/\alpha}(m_1\wedge m_2))$
and $q=\dfrac{\log (m)}{\log (m)-1}$. By H{\"o}lder's inequality we get
\begin{equation}\label{estEnormMrank}
\begin{split}
\bE\left (\Delta^{2}W\right)\leq \left (\bE \Delta^{2\log m}\right )^{1/\log m} \left (\bE\,W^{q}\right )^{1/q}.
\end{split}
\end{equation}
We first estimate $\left (\bE \Delta^{2\log m}\right )^{1/\log m}$. Inequalities \eqref{proba_1} and \eqref{proba_2} imply that 
\begin{equation}\label{estEM}
\begin{split}
\left (\bE \Delta^{2\log m}\right )^{1/\log m} &=\left ( \overset{+\infty}{\underset{0}{\int}}\mathbb P \left (\Delta> t^{1/(2\log m)}\right )\mathrm{d}t\right )^{1/\log m}\\&\leq
\left (m \overset{(t^{*})^{2k}}{\underset{0}{\int}}\exp\{-t^{1/\log m}\nu_1\}\mathrm{d}t\right .\\&\hskip 0.5 cm+\left .3\sqrt{m} \overset{+\infty}{\underset{(t^{*})^{2k}}{\int}}\exp\{-t^{1/(2\log m)}\nu_2\}\mathrm{d}t
\right )^{1/\log m}\\&\hskip -1.5 cm\leq e\left (\log (m)\nu_1^{-\log m}\Gamma(\log m)+\dfrac{6}{\sqrt{m}}\log(m)\; \nu_2^{-2\log m}\Gamma(2\log m)\right )^{1/\log m}.
%\\&=\left ( \overset{+\infty}{\underset{0}{\int}} 2\log(m)\; m\; t^{2\log( m)-1} e^{-t^{2}\nu^{2}} \mathrm{d}t\right )^{1/\log m
\end{split}
\end{equation}
The Gamma-function satisfies the following bound: 
\begin{equation}\label{Gamma}
\text{for}\quad x\geq 2,\quad\Gamma(x)\leq \left (\dfrac{x}{2}\right )^{x-1}.
\end{equation} 
We give a proof of this inequality in the Appendix. Plugging it into \eqref{estEM} we compute  
\begin{equation*}
\begin{split}
\left (\bE \Delta^{2\log m}\right )^{1/\log m}& \leq e\left ((\log (m))^{\log m}\nu_1^{-\log m}2^{1-\log m}\right .\\&+\left .\dfrac{6}{\sqrt{m}}(\log (m))^{2\log m}\nu_2^{-2\log m}\right )^{1/\log m}. %\exp\left \{1-\dfrac{\log 2}{2}\right \}.
\end{split}
\end{equation*}
Observe that $n>n^{*}$ implies 
$\nu_1\log m\leq\nu_2^{2}$ and we obtain
 \begin{equation}\label{estEM-1}
 \begin{split}
 \left (\bE \Delta^{2\log m}\right )^{1/\log m}& \leq e\log (m)\nu_1^{-1}\left (2^{1-\log m}+\dfrac{6}{\sqrt{m}}\right )^{1/\log m}.
 \end{split}
 \end{equation}
If $\bE(W^{q})<1$ we get \eqref{estEnormMrank-1} directly from \eqref{estEnormMrank}.
If $\bE(W^{q})\geq 1$, the bound $W\leq m_1\wedge m_2$ implies that
\begin{equation}\label{rankB}
\begin{split}
\left (W\right )^{1/(\log (m)-1)}&\leq \exp \left \{\dfrac{\log(m_1\wedge m_2)}{\log (m)-1}\right \}
\leq \exp \left \{\dfrac{\log(m)-\log 2}{\log (m)-1}\right \} \\&\leq e\exp\left \{\dfrac{1-\log2}{\log (m)-1}\right \}. %\exp \left \{1+\dfrac{\log 2}{2}\right \}.
\end{split}
\end{equation}
and we compute
\begin{equation}\label{estErank^q}
\begin{split}
\left (\bE\left (W^{1+1/(\log (m)-1)}\right )\right )^{1-1/\log m}&\leq \bE\left (W^{1+1/(\log (m)-1)}\right )\\&\leq e\exp\left \{\dfrac{1-\log2}{\log (m)-1}\right \}%\exp \left \{1+\dfrac{\log 2}{2}\right \}
\bE( W).
\end{split}\end{equation}
The function
\begin{equation*}
\begin{split}
\left (2^{1-\log m}+\dfrac{6}{\sqrt{m}}\right )\exp\left \{\dfrac{(1-\log2)\log m}{\log (m)-1}\right \}\\
=\dfrac{e}{2}\left (2^{1-\log m}+\dfrac{6}{\sqrt{m}}\right )\exp\left \{\dfrac{1-\log2}{\log (m)-1}\right \}
\end{split}
\end{equation*}
is a decreasing function of $\log m$ which is smaller then $1$ for $\log m\geq 5$. This implies
\begin{equation}\label{constfinal}
\left (2^{1-\log m}+\dfrac{6}{\sqrt{m}}\right )^{1/\log m}\exp\left \{\dfrac{1-\log2}{\log (m)-1}\right \}<1
\end{equation}
Plugging \eqref{estErank^q} and  \eqref{estEM-1} into \eqref{estEnormMrank} and using \eqref{constfinal} 
 we get \eqref{estEnormMrank-1}.
This completes the proof of Lemma \ref{lemma2}.
\end{proof}

%The condition ``for any $t>0$ we have $\Delta\leq t$ with  a probability at least $1-me^{-t^{2}(m_1\wedge m_2)n}$'' is particularly satisfied in the case of matrix completion (see Lemma \ref{lemma}). 

 The natural choice of $t$ in  Lemma \ref{lemma1} is of the order $\log m$ (see the discussion in \cite{Koltchinskii-Tsybakov}).
Then, in Theorems \ref{thm1} and \ref{thm2} we can take  $\sqrt{\lambda}=2\varrho c\sqrt{\dfrac{(m_1\vee m_2)\log (m)}{n}}$, where the constant $c$ is large enough, to obtain the following corollary.
%For $n>(m_1\wedge m_2)\log ^{1+2/\alpha}(m)$ in the case of sub-exponential noise and  matrices with uniformly bounded entries (resp. $n>4(m_1\wedge m_2)\log(m)$ in the case of statistical learning) we can take in the Theorem \ref{thm1} $\sqrt{\lambda}=2kc\sqrt{\dfrac{(m_1\vee m_2)\log (m)}{n}}$, where $c$ equals either $\sigma\wedge a$ or $\eta$, to obtain the following corollary:
% % % % % % % % % % % % % % % % % % % % % % % % % % % % % % % % % % % % % % % % % % % % %
\begin{Corollary}\label{cor1} Let one of the set of conditions (\textbf A) or (\textbf B) be satisfied. Assume $n>n^{*}$, $\log m\geq 5$, $\varrho\geq 1$ and $\sqrt{\lambda}=2\varrho c\sqrt{\dfrac{(m_1\vee m_2)\log (m)}{n}}$. Then, 
\begin{enumerate}
\item[(i)] with probability at least $1-3/(m_1+m_2)$, one has
\begin{equation*}
\begin{split}
\dfrac{\parallel \hat{A}-A_{0}\parallel _{2}}{\sqrt{m_1m_2}}&\leq \underset{A\in \mathbb{R}^{m_{1}\times m_{2}}}{\mathrm{inf}}\left  \{\dfrac{ \parallel A-A_{0}\parallel _{2}}{\sqrt{m_1m_2}}+2\sqrt{\lambda\max\left (\dfrac{\rank (A)_{0}}{\varrho^{2}},\rank (A)\right )}\right  \}
 \end{split}
\end{equation*}
and, in particular,
\begin{equation*} %\label{concn}
\begin{split}
\dfrac{\parallel \hat{A}-A_{0}\parallel _{2}}{\sqrt{m_1m_2}}\leq 4c\sqrt{\dfrac{(m_1\vee m_2)\log (m)\rank (A_0)}{n}},
 \end{split}
\end{equation*}
\item  [(ii)] with probability at least $1-3/(m_1+m_2)$, one has
\begin{equation*}
\begin{split}
\dfrac{\parallel \hat A -A_0\parallel _{2}^{2}}{m_1m_2}\leq \underset{A\in \mathbb{R}^{m_{1}\times m_{2}}}{\mathrm{inf}}
\left \{\left (\dfrac{2\varrho^{2}+1}{2\varrho^{2}-1}\right )\dfrac{ \parallel A-A_0 \parallel_{2}^{2} }{m_1m_2}+\dfrac{4\varrho^{2}\lambda}{2\varrho^{2}-1} \rank (A)\right \},
\end{split}
\end{equation*}
\item[(iii)] 
\begin{equation*}
\begin{split}
\dfrac{\bE\parallel \hat{A}-A_{0}\parallel _{2}}{\sqrt{m_1m_2}}&\leq \underset{A\in \mathbb{R}^{m_{1}\times m_{2}}}{\mathrm{inf}}\left  \{\dfrac{\parallel A-A_{0}\parallel _{2}}{\sqrt{m_1m_2}}\right .\\&\hskip 0.5 cm\left .+\dfrac{5}{2}\sqrt{\lambda\max\left (\rank (A),\dfrac{\rank (A)_{0}}{\varrho^{2}},\dfrac{1}{4\varrho^{2}}\right )}\right  \}
 \end{split}
\end{equation*}
and, in particular,
\begin{equation*}
\begin{split}
\dfrac{\bE\parallel \hat{A}-A_{0}\parallel _{2}}{\sqrt{m_1m_2}}\leq 5c\sqrt{\dfrac{(m_1\vee m_2)\log (m)}{n}\max\left (\rank (A_0),\dfrac{1}{4}\right)},
 \end{split}
\end{equation*}
\item  [(iv)] 
\begin{equation*}
\begin{split}
\dfrac{\bE\parallel \hat A -A_0\parallel _{2}^{2}}{m_1m_2}&\leq \underset{A\in \mathbb{R}^{m_{1}\times m_{2}}}{\mathrm{inf}}
\Big\{\left (\dfrac{2\varrho^{2}+1}{2\varrho^{2}-1}\right )\dfrac{ \parallel A-A_0 \parallel_{2}^{2} }{m_1m_2}\\
&\hskip 0.5 cm+\left (\dfrac{4\varrho^{2}\lambda}{2\varrho^{2}-1}\right )\max\left (\rank (A),\dfrac{1}{2}\right )\Big\},
\end{split}
\end{equation*}
and, in particular,
\begin{equation*}
\begin{split}
\dfrac{\bE\parallel \hat{A}-A_{0}\parallel _{2}^{2}}{m_1m_2}\leq \dfrac{16 \,c^{2}(m_1\vee m_2)\log (m)}{n}\max\left (\rank (A_0),\dfrac{1}{2}\right ).
 \end{split}
\end{equation*}
\item[(v)]  with probability at least $1-3/(m_1+m_2)$, one has
\begin{equation*}
 \left \| \hat A -A_0\right \|_{\infty}\leq 3\rho c \sqrt{m_1m_2\dfrac{(m_1\wedge m_2)\log m}{n}}
\end{equation*}
\item[(vi)] with probability at least $1-3/(m_1+m_2)$, one has
\begin{equation*}
\begin{split}
\dfrac{\parallel \hat A -A_0\parallel _{2}^{2}}{m_1m_2}&\leq\left (\dfrac{2\varrho^{2}+1}{2\varrho^{2}-1}\right ) \underset{0<q\leq 2}{\inf}\dfrac{\lambda^{1-q/2}\parallel A_0 \parallel_{q}^{q}}{(m_1m_2)^{q/2}},
\end{split}
\end{equation*}

\end{enumerate}

\end{Corollary}
\begin{proof}
(i) - (iv) are straightforward in view of Theorems \ref{thm1} and \ref{thm2}. (v) is a consequence of Theorem \ref{thm3} (iii). The proof of (vi) follows from (i) using the same argument as in \cite{Koltchinskii-Tsybakov} Corollary 2. 
\end{proof}
This corollary guarantees that the normalized Frobenius error $\dfrac{\parallel \hat{A}-A_{0}\parallel _{2}}{\sqrt{m_1m_2}}$ of the estimator $\hat A$ is small whenever $n>C(m_1\vee m_2)\log(m)\rank (A_0)$ with a constant $C$ large enough. This quantifies the sample size $n$ necessary for successful matrix completion from noisy data. 

 Comparing Corollary \ref{cor1} with Theorem 6 and Theorem 7 of \cite
 {Koltchinskii-Tsybakov} we see that, in the case of Gaussian errors and for the statistical learning setting, the rate of convergence of our estimator is optimal,  for the class of matrices $\mathcal{A}(r,a)$ defined as follows: for any $A_0\in \mathcal{A}(r,a)$ the rank of $A_0$ is supposed not to be larger than a given $r$ and all the entries of $A_0$ are supposed to be bounded in absolute value by a constant $a$.
 
   \section{Matrix Regression}\label{matrix regression}
   In this section we apply the general oracles inequalities of Theorems \ref{thm1} and \ref{thm2} to the matrix regression model and compare our bounds to those obtained by Bunea, She and Wegkamp in \cite {2010arXiv1004.2995B}.  Recall that matrix regression model is given by 
   \begin{equation}\label{regression_bis}
   U_i=V_i\,A_0+E_i\qquad i=1,\dots, l,
   \end{equation}
   where $U_i$ are $1\times m_2$ vectors of response variables, $V_i$ are $1\times m_1$ vectors of predictors, $A_0$ is a unknown $m_1\times m_2$ matrix of regression coefficients and $E_i$ are random $1\times m_2$ vectors of noise with independent entries and mean zero.
   
    As mentioned in the section \ref{def}, we can equivalently write this model as a trace regression model. Let $U_i=(U_{ik})_{k=1,\dots, m_2},E_i=(E_{ik})_{k=1,\dots, m_2}$ and $Z^{T}_{ik}=e_k(m_2)\,V_i$ where $e_k(m_2)$ are the $m_2\times 1$ vectors of the canonical basis of $\mathbb {R}^{m_2}$. Then we can write \eqref{regression_bis} as
   \begin{equation*}
   U_{ik}=\mathrm{tr}(Z_{ik}^{T}A_0)+E_{ik}\qquad i=1,\dots, l \quad\text{and}\quad k=1,\dots, m_2.
   \end{equation*}
   Set 
   %$n=lm_2$ and 
   $V=\left  (V_1^{T},\dots, V_l^{T}\right  )^{T}$, $U=\left  (U_1^{T},\dots, U_l^{T}\right  )^{T}$ and $E=\left  (E_1^{T},\dots, E_l^{T}\right  )^{T}$, then for deterministic predictors
   \begin{equation*}
   \parallel A \parallel _{L_{2}(\Pi)}^{2}=\dfrac{1}{l\,m_2}\parallel VA\parallel^{2}_2
   \end{equation*}
   Note that we use  Assumption \ref{ass1} in the proof of Theorem \ref{thm1} to derive \eqref{triangle} from \eqref{estim}. In the case of matrix regression with deterministic $V_i$, we do not need this assumption and proceed as follows. Let $\mathcal{P}_{A}$ denote the orthogonal projector on the linear span of the columns of matrix $A$ and  let $\mathcal{P}_{A}^{\bot}=1-\mathcal{P}_{A}$. Note that $A\mathcal{P}_{A^{T}}^{\bot}=0$. Then,  one has $M=V^{T}E=V^{T}\mathcal{P}_{V}E$. Now, we use \eqref{estim} and the fact that
   \begin{equation*}
      \begin{split}
      \langle M,A-\hat A\rangle=\langle V^{T}\mathcal{P}_{V}E,A-\hat A\rangle=\langle \mathcal{P}_{V}E,V(A-\hat A)\rangle.
      \end{split}
      \end{equation*}
      Hence, the trace duality yields \eqref{triangle} where we set $\Delta=\left \| \mathcal{P}_{V}E\right\|_\infty$. Thus, in the case of matrix regression with deterministic $V_i$, we have proved that Theorems \ref{thm1} and \ref{thm2} hold with $\Delta=\left \| \mathcal{P}_{V}E\right\|_\infty$ even if Assumption \ref{ass1} is not satisfied.
%      Which implies that in the case of matrix regression with deterministic predictors, we  instead of Assumption \ref{ass1},  it is enough to have that
%      \begin{equation}\label{ass-bis}
%            \parallel \mathcal{P}_{V^{T}}A\parallel_{2}\leq \mu \parallel VA\parallel_{2}.
%            \end{equation}
%            This condition is satisfied for any $A$ if $\mu$ is larger then the inverse of the smallest positive singular value of $V$ i.e.
%            \begin{equation*}
%            \mu \geq \sigma_{\rank (V)}^{-1}(V).
%            \end{equation*}
%         Using that $M=V^{T}\mathcal{P}_{V}E$ we get the following upper bound on $ \Delta$ :
%         \begin{equation}
%         \Delta\leq \sigma_1(V) \left \| \mathcal{P}_{V}E\right\|_\infty.
%         \end{equation}

         In order to get an upper bound on $ S=\underset{W\in\mathcal B_{m_1\wedge m_2}}{\sup}\dfrac{\bE\left (\Delta^{2}W\right)}{ \max\{\bE(W), 1\}}$ in the case of Gaussian noise we will use the following result. 
 \begin  {lemma}[\cite {2010arXiv1004.2995B}, Lemma 3]\label{lemma_bunea}
         Let $r=\rank(V)$ and assume that $E_{ij}$ are independent $N(0,\sigma^{2})$ random variables.Then
           \begin{equation}\label{lemma_bunea_E}
          \bE (\left \| \mathcal{P}_{V}E\right\|_\infty)\leq \sigma(\sqrt{m_2}+\sqrt{r})
          \end{equation}
          and\begin{equation}\label{lemma_bunea_P}
          \mathbb P\left \{\left \| \mathcal{P}_{V}E\right\|_\infty\geq \bE (\left \| \mathcal{P}_{V}E\right\|_\infty)+\sigma t\right \}\leq \exp\left \{-t^{2}/2\right \};
          \end{equation}
\end{lemma} 
          Using \eqref{lemma_bunea_P} and Lemma \ref{concentratio_bound} applied to $\mathcal{P}_{V}E$ we get the following bound on $S$:
          \begin{lemma}
          Assume that $E_{ij}$ are independent $N(0,\sigma^{2})$, then\begin{equation}
          S\leq 2\left (\bE \left (\left \| \mathcal{P}_{V}E\right\|_\infty\right )\right)^{2}+4\,e\,\sigma^{2}(2\log(m_1\wedge m_2)+1).
          \end{equation}
          \end{lemma}
          %Let $C_{0}(V)=\sigma_1(V)/\sigma_{\rank(V)}(V)$. 
          For $\log m_2\geq 4$, we have that $m_2\geq 2\,e(2\log m_2+1)$. Then, these two lemmas imply that in Theorems \ref{thm1} and \ref{thm2} we can take  $\sqrt{\lambda}=4\sigma \left (\sqrt{r}+\sqrt{m_2}\right )$ to get the following corollary:
          \begin{Corollary}
          Assume that $E_{ij}$ are independent $N(0,\sigma^{2})$, $\log m_2\geq 4$, $\rho\geq 1$ and $\sqrt{\lambda}=4\rho\sigma \left (\sqrt{r}+\sqrt{m_2}\right )$. Then
          \begin{enumerate}
          \item[(i)] with probability at least $1-\exp\left (-\dfrac{m_2+r}{2}\right )$, one has
          \begin{equation*}
          \parallel V(\hat{A}-A_{0})\parallel _{2} \leq \underset{A\in \mathbb{R}^{m_{1}\times m_{2}}}{\mathrm{inf}}\left  \{\parallel V\left (A-A_{0}\right )\parallel _{2}+2\sqrt{\lambda\max\left (\dfrac{\rank (A)_{0}}{\varrho^{2}},\rank (A)\right )}\right  \}
          \end{equation*}
          and, in particular,
        \begin{equation*} %\label{concn}
          \parallel V(\hat{A}-A_{0} ) \parallel _{2} \lesssim\left( \sqrt{r}+\sqrt{m_2}\right) \sqrt{\rank (A_0)},
        \end{equation*}
          \item  [(ii)] with probability at least $1-\exp\left (-\dfrac{m_2+r}{2}\right )$, one has
          \begin{equation*}
          \begin{split}
          \parallel V(\hat{A}-A_{0})\parallel _{2}^{2}\leq \underset{A\in \mathbb{R}^{m_{1}\times m_{2}}}{\mathrm{inf}}
          \left \{\left (\dfrac{2\varrho^{2}+1}{2\varrho^{2}-1}\right )\parallel V(A-A_{0})\parallel _{2}^{2}+\dfrac{4\varrho^{2}\lambda}{2\varrho^{2}-1} \rank (A)\right \},
          \end{split}
          \end{equation*}
          \item[(iii)] 
          \begin{equation*}
          \begin{split}
          \bE\left (\parallel V(\hat{A}-A_{0})\parallel _{2}\right )&\leq \underset{A\in \mathbb{R}^{m_{1}\times m_{2}}}{\mathrm{inf}}\Big  \{\parallel V(A-A_{0})\parallel _{2}\Big.\\&\hskip 0.5 cm\left .+\dfrac{5}{2}\sqrt{\lambda\max\left (\rank (A),\dfrac{\rank (A)_{0}}{\varrho^{2}},\dfrac{1}{4\varrho^{2}}\right )}\right  \}
           \end{split}
          \end{equation*}
          and, in particular,
          \begin{equation*}
          \begin{split}
          \bE\left (\parallel V(\hat{A}-A_{0})\parallel _{2}\right )\lesssim\left( \sqrt{r}+\sqrt{m_2}\right) \sqrt{\max\left (\rank (A_0),1/4\right )}
           \end{split}
          \end{equation*}
          \item  [(iv)] 
          \begin{equation*}
          \begin{split}
          \bE\left (\parallel V(\hat{A}-A_{0})\parallel _{2}^{2}\right )&\leq \underset{A\in \mathbb{R}^{m_{1}\times m_{2}}}{\mathrm{inf}}
          \Big\{\left (\dfrac{2\varrho^{2}+1}{2\varrho^{2}-1}\right )\parallel V(A-A_{0})\parallel _{2}^{2}\\
          &\hskip 0.5 cm+\left (\dfrac{4\varrho^{2}\lambda}{2\varrho^{2}-1}\right )\max\left (\rank (A),\dfrac{1}{2}\right )\Big\},
          \end{split}
          \end{equation*}
          and, in particular,
          \begin{equation*}
                    \bE\left (\parallel V(\hat{A}-A_{0})\parallel _{2}^{2}\right )\lesssim\left( \sqrt{r}+\sqrt{m_2}\right)^{2} \max\left (\rank (A_0),1/2\right )
                    \end{equation*}
          \end{enumerate}
          The symbol $\lesssim$ means that the inequality holds up to multiplicative numerical constants.
         \end{Corollary}
         This Corollary shows that our error bounds are comparable to those obtained in \cite {2010arXiv1004.2995B}. Points (i) and (iii) are new; here we have inequalities with leading constant 1. The results (ii) and (iv) give the same bounds as in \cite {2010arXiv1004.2995B} up to constants and to an  additional exponentially small term in the analog of (iv) in \cite {2010arXiv1004.2995B}.

  \section{Appendix}
   For completeness, we give here the proof of \eqref{Gamma}.
   \begin{proposition}\label{proposition1}
    $$\Gamma(x)\leq \left (\dfrac{x}{2}\right )^{x-1}\;\text{for}\;x\geq 2$$
   \end{proposition}
   \begin{proof}
   We set $\tilde{\Gamma}(x)=\Gamma(x)\left (\dfrac{2}{x}\right )^{x-1}.$ Using functional equation for $\Gamma$ we note that $\tilde{\Gamma}(x)=\tilde{\Gamma}(x+1)2^{-1}\left (1+\dfrac{1}{x}\right )^{x}.$ Applying this equality $n$ times we get
   \begin{equation}\label{gamma-tilde}
   \tilde{\Gamma}(x)=\tilde{\Gamma}(x+n)\exp\left \{\underset{j=0}{\overset{n-1}{\Sigma}}(x+j)\log \left (\dfrac{x+j+1}{x+j}\right )-n\log 2\right \}.
   \end{equation}
   By Stirling's formula, we have that 
       \begin{equation*} %\label{stirling}
       \tilde{\Gamma}(x+n)=\sqrt{\dfrac{\pi}{2}}\left (\dfrac{2}{e}\right )^{x+n}\sqrt{x+n}\left (1+O\left (\dfrac{1}{x+n}\right )\right ).
       \end{equation*}
       Plugging this into \eqref{gamma-tilde} we obtain
       \begin{equation*}
       \begin{split}
       \log \tilde{\Gamma}(x)&= \underset{n\rightarrow \infty}{\lim}\left [\underset{j=0}{\overset{n-1}{\Sigma}}\left ((x+j)\log \left (\dfrac{x+j+1}{x+j}\right )\right )-n\right . \\&\left.-n\log 2+\dfrac{1}{2}\log(x+n)\right] +x(\log 2-1)+\dfrac{1}{2}\log(\pi/2).
       \end{split}
       \end{equation*}
   Note that 
   \begin{equation*} %\label{telescop}
   \begin{split}
   \underset{j=0}{\overset{n-1}{\Sigma}}\left ((x+j)\log \left (\dfrac{x+j+1}{x+j}\right )\right )-n&=\underset{j=0}{\overset{n-1}{\Sigma}}\left (\underset{0}{\overset{1}{\int}}\dfrac{x+j}{x+j+u}\,\mathnormal{du}-1\right )\\&=-\underset{j=0}{\overset{n-1}{\Sigma}}\underset{0}{\overset{1}{\int}}\dfrac{u}{x+j+u}\,\mathnormal{du}.
   \end{split}
   \end{equation*}

    Defining 
   \begin{equation*}
   \begin{split}
   F(x)&=\underset{n\rightarrow \infty}{\lim}\left (-\underset{j=0}{\overset{n-1}{\Sigma}}\underset{0}{\overset{1}{\int}}\dfrac{u}{x+j+u}\,\mathnormal{du} +\dfrac{1}{2}\log(x+n)-n\log2\right )
   \end{split}\end{equation*}
   we have
   \begin{equation}\label{gamma-F}
   \log \tilde{\Gamma}(x)=F(x)+x(\log2 -1)+\dfrac{1}{2}\log\left (\dfrac{\pi}{2}\right ).
   \end{equation}
   Observe that $F$ is infinitely differentiable on $[1,+\infty)$. Moreover the series defining $F$ can be differentiated $k$ times to obtain $F^{(k)}$. Thus 
   \begin{equation}\label{Fprime}
   \begin{split}
   F'(x)&=\underset{n\rightarrow \infty}{\lim}\left (\underset{j=0}{\overset{n-1}{\Sigma}}\underset{0}{\overset{1}{\int}}\dfrac{u}{(x+j+u)^{2}}\,\mathnormal{du}+\dfrac{1}{2(x+n)}\right )\\
   &=\underset{n\rightarrow \infty}{\lim}\underset{j=0}{\overset{n-1}{\Sigma}}\underset{0}{\overset{1}{\int}}\dfrac{u}{(x+j+u)^{2}}\,\mathnormal{du}
   \end{split}
   \end{equation} 
   and \begin{equation}
   F''(x)=-\underset{n\rightarrow \infty}{\lim}\underset{j=0}{\overset{n-1}{\Sigma}}\underset{0}{\overset{1}{\int}}\dfrac{2u}{(x+j+u)^{3}}\,\mathnormal{du}<0.
   \end{equation}
    The relation \eqref{gamma-F} implies that $(\log \tilde{\Gamma})'(2)=F'(2)+\log2-1$. Using \eqref{Fprime} for $x=2$ we get
   \begin{equation*}
   \begin{split}
   (\log \tilde{\Gamma})'(2)&=\underset{n\rightarrow \infty}{\lim}\left (\log(n+2)-\log2 -\underset{j=0}{\overset{n-1}{\Sigma}}\dfrac{1}{j+3}\right )+\log2-1\\
      &=\underset{n\rightarrow \infty}{\lim}\left (\log n-\underset{j=1}{\overset{n}{\Sigma}}\dfrac{1}{j} \right )+\dfrac{1}{2}=-\gamma+\dfrac{1}{2}<0
   \end{split}\end{equation*}
   where $\gamma$ is the Euler's constant. Together with $\log \tilde{\Gamma}(2)=0$ and $(\log \tilde{\Gamma})''(2)=F''(2)<0$ this implies that
   \begin{equation*}
   \log \tilde{\Gamma}(x)<0
   \end{equation*}
   for any $x\geq 2$. This completes the proof of Proposition \ref{proposition1}.
   \end{proof}   % % % % % % % % % % % % % % % % % % % % % % % % % % % % % % % % % % % % % % % % % % % % % % % % % % % % % % % % % %

\textbf{Acknowledgements}.
 It is a pleasure to thank A. Tsybakov for introducing me to these interesting subjects and for guiding this work.
%\bibliographystyle{plain}
%\bibliography{biblio} 

\end{document}